\documentclass [11pt,a4paper]{article}

\usepackage{amsmath}
\usepackage{amsthm}
\usepackage{graphics}
\usepackage{latexsym}
\usepackage{amssymb}
\usepackage{enumerate}
\usepackage[latin2]{inputenc}
\usepackage{color}

\newtheorem  {lemma}         {Lemma}

\newtheorem* {theorem*}      {Theorem}
\newtheorem* {lemma*}        {Lemma}
\newtheorem* {corollary*}    {Corollary}
\newtheorem* {proposition*}  {Proposition}
\newtheorem* {definition*}   {Definition}
\newtheorem* {remark*}       {Remark}
\newtheorem* {remarks*}      {Remarks}

\newcounter{aux}

\def \N {\mathbb N}
\def \Z {\mathbb Z}

\def \R {\mathbb R}

\def \T {\mathbb T}

\def \Ordo {{\cal O}}
\def \ind{1\!\!1}

\def \dsum {\displaystyle\sum}

\def \del   {\delta}

\def \veps  {\varepsilon}

\def \sig    {\sigma}

\def \ome   {\omega}
\def \om {\omega}

\def \Ups   {\Upsilon}

\def \Om   {\Omega}

\def\uome{\underline{\omega}}
\def\uo{\uome}

\def\wt{\widetilde}

\def\old{\overline{\delta}}
\def\ols{\overline{\sigma}}

\def\beqs{\begin{eqnarray*}}
\def\eeqs{\end{eqnarray*}}
\def\beq{\begin{eqnarray}}
\def\eeq{\end{eqnarray}}
\def\beas{\begin{eqnarray*}}
\def\eeas{\end{eqnarray*}}
\def\bea{\begin{eqnarray}}
\def\eea{\end{eqnarray}}

\def \prob        {\ensuremath{\mathbf{P}}}
\def \expect      {\ensuremath{\mathbf{E}}}
\def \var         {\ensuremath{\mathbf{Var}}}
\def \cov         {\ensuremath{\mathbf{Cov}}}

\def\pt{\partial_t}
\def\px{\partial_x}
\newcommand{\abs}[1]{\left|{#1}\right|}

\def \b {\beta}

\def \Omn {\Omega^{n}}

\def \Tn {\T^{n}}
\def \Ln {L^{n}}

\def \grn {\nabla^{n}}
\def \gradn {\grn}

\def \Xn {{\cal X}^{n}}

\def \mun  {\mu^{n}}
\def \munt {\mu^{n}_t}

\def \nun  {\nu^{n}}

\def \pin  {\pi^{n}}

\def\taun {\tau^{n}}

\def \g {\gamma}
\def \DD {\mathcal{D}}

\def \nuun {{\widetilde{\nu}}^{n}}
\def \nuunt {{\widetilde{\nu}}^{n}_t}
\def \tetn {{\widetilde{\theta}}^{(n)}}

\def \taun {{\widetilde{\tau}}^{(n)}}
\def \uun  {{\widetilde{u}}^{(n)}}
\def \vvn  {{\widetilde{v}}^{(n)}}
\def \ffnt  {\widetilde{f}^{n}_t}
\def \ffn  {\widetilde{f}^{n}}
\def \Err   {\textup{Err}}
\def \RR {\textup{R}}

\newcommand{\mb}[1] {\mathbf{#1}}

\def \wzeta {\widehat{\zeta}^{l}}
\title
{Hydrodynamic limit for perturbation of a hyperbolic equilibrium
point in two-component systems}

\author {Benedek Valkó \footnote{\textsc {\scriptsize{Alfréd Rényi Institute of Mathematics, Hungarian Academy of
Sciences, Reáltanoda u.~13-15, H-1053, Budapest, Hungary and
Institute of Mathematics, Technical University Budapest, Egry
J\'ozsef u.~1, H-1111 Budapest, Hungary. Email:{\tt
valko{@}renyi.hu}}}}}

\begin{document}

\setlength{\baselineskip}{1.23\baselineskip}

\maketitle
\begin{abstract}
We consider one-dimensional, locally finite interacting particle
systems with two conservation laws. The models have a family of
stationary measures with product structure and we assume the
existence of a uniform bound on the inverse of the spectral gap
which is quadratic in the size of the system. Under Eulerian
scaling the hydrodynamic limit for the macroscopic density
profiles leads to a two-component system of conservation laws. The
resulting pde is hyperbolic inside the physical domain of the
macroscopic densities, with possible loss of hyperbolicity at the
boundary.

We investigate the propagation of small perturbations around a
\emph{hyperbolic} equilibrium point. We prove that the
perturbations essentially evolve according to two \emph{decoupled}
Burgers equations. The scaling is not Eulerian: if the lattice
constant is $n^{-1}$, the perturbations are of order $n^{-\beta}$
then time is speeded up by $n^{1+\b}$. Our derivation holds for
$0<\beta< \frac15$. The proof relies on Yau's relative entropy
method, thus it applies only in the regime of smooth solutions.

This result is an extension of \cite{seppalainen} and
\cite{tothvalko1} where the analogue result was proved for systems
with one conservation law. It also complements \cite{tothvalko3}
where it was shown that perturbations around a non-hyperbolic
boundary equilibrium point are driven by a universal two-by-two
system of conservation laws.

\end{abstract}
%%%%%%%%%%%%%%%%%%%%%%%%%%%%%%%%%%%%%%%%%%%%%%%%%%%%%%%%%%%%%%%%
%%%%%%%%%%%%%%%%%%%%%%%%%%%%%%%%%%%%%%%%%%%%%%%%%%%%%%%%%%%%%%%%
\section{Introduction}
\label{section:hyp_intro}

There are several results dealing with the perturbation analysis
of hydrodynamic limits for interacting particle systems. In the
landmark paper \cite{espositomarrayau} the authors prove that for
the asymmetric simple exclusion, in dimensions higher than 2,
perturbations of order $n^{-1}$ of a constant profile evolve
according to a certain parabolic equation under diffusive scaling
(time rescaled by $n^2$, space by $n$). It is well-known, that
under Eulerian scaling (time rescaled by $n$, space by $n$) the
hydrodynamic limit leads to a hyperbolic conservation law (the
Burgers equation), the perturbation limit gives the same equation
with the Navier-Stokes correction. (For a survey on the
microscopic interpretations of the Navier-Stokes equations see the
end of Chapter 7 of \cite{kipnislandim}.)

Motivated by \cite{espositomarrayau} T.~Sepp\"al\"ainen
investigated a similar problem in one dimension for the so-called
totally asymmetric stick process. In \cite{seppalainen} he proves
that an $\Ordo(n^{-\beta})$ perturbation of the constant profile
is governed by the Burgers equation (even after the appearance of
shocks) if time is rescaled by $n^{1+\beta}$ and space by $n$,
where $\beta\in (0,\frac12)$ is a fixed constant. Independently,
in \cite{tothvalko1} the authors partially extend this result by
proving that one gets \emph{universally} the Burgers equation in
the hydrodynamic limit for similar perturbations of equilibrium
for a wide class of one-dimensional interacting particle systems
with one conservation law. The models are not reversible and not
necessarily attractive. The proof relies on H.~T.~Yau's relative
entropy method, it only applies in the smooth regime of solutions
and it only works for $\b\in (0,\frac15)$. It is conjectured that
the result should hold for all $\b\in (0,\frac12)$ even without
the smoothness condition as in the result of \cite{seppalainen}.

This universal result may be explained by the following arguments.
Under Eulerian scaling these systems admit in the hydrodynamic
limit a hyperbolic conservation law  of the form
\beq%
\pt u+\px J(u)=0. \label{eq:hyp_pde_euler_1}
\eeq%
Taking a point $u_0$ with $J''(u_0)\neq 0$ simple (although
formal) calculations yield that solutions of
(\ref{eq:hyp_pde_euler_1}), with initial conditions which are
small perturbations $u_0$, are governed by the Burgers equation.
See \cite{tothvalko1} for the 'more precise' formulation.

In the present paper we give an extension of the results of
\cite{seppalainen,tothvalko1} for systems with 2 conserved
quantities. In \cite{tothvalko2} a general one-dimensional family
of lattice-models was introduced. The models are locally finite
interacting particle systems with two conservation laws which
possess a family of stationary measures with product structure. In
that paper it is shown (in the regime of smooth solutions)
that in Eulerian scaling we get a hydrodynamic limit of the form%
\beq%
\left\{\begin{array}{rcl}
  \pt u+\px \Phi(u,v)&=&0,
\\[3pt]
  \pt v+\px \Psi(u,v)&=&0,
\\
\end{array}
\right.\label{eq:hyp_pde_euler_intro}
\eeq%
where $(u,v)\in \DD$ and $\DD$  is a convex compact polygon, the
the physical domain (see (\ref{eq:hyp_DD}) for the definition). We
also note, that \cite{ollavaradhanyau} gives the first major
result about the Eulerian hydrodynamic limit of multi-component
hyperbolic systems. In \cite{tothvalko2} it was also shown that an
Onsager-type symmetry relation holds for the macroscopic flux
functions $\Phi, \Psi$ (see Lemma \ref{lemma:hyp_onsager}). One of
the consequences of this relation is that inside the physical
domain $\DD$ the pde is (weakly) hyperbolic, i.e. the Jacobian can
be diagonalized in the real sense. Experience shows, that the
limiting pde is strongly hyperbolic (the Jacobian has two distinct
real eigenvalues) in the whole physical domain except some special
points on the boundary $\partial \DD$.

We consider perturbations of order $n^{-\b}$ around a constant
equilibrium point $(u_0,v_0)\in \DD$, which is \emph{strictly}
hyperbolic. We prove that rescaling time by $n^{1+\b}$ and space
by $n$ the evolution of the perturbations are governed by two
\emph{decoupled} equations. (These are 'usually' Burgers
equations, see the remark at the end of subsection
\ref{subs:hyp_formal_pert}.) This result agrees with the formal
perturbation of the pde (\ref{eq:hyp_pde_euler_intro}) e.g.~with
the method of weakly nonlinear geometric optics (see
\cite{dipernamajda,hunterkeller}).

The reason for the decoupling of the resulting pde system is the
strict hyperbolicity, basically, the two different eigenvalues
(sound speeds) cause the equations to separate.  In the paper
\cite{tothvalko3} perturbation around a special
\emph{non-hyperbolic} point was considered in a similar setting,
it was proved that in that case in the limit the evolution obeys a
two-by-two system of conservation laws which cannot be decoupled.
The treatment of that problem needs more complex tools than our
proofs, sophisticated pde methods are used besides Yau's method.

Our proof follows the relative entropy method using similar steps
as \cite{tothvalko1} (thus it only applies in the regime of smooth
solutions), but it also heavily relies on the Onsager-type
symmetry relation proved in \cite{tothvalko2}. We assume the
existence of a uniform  bound on the inverse of the spectral gap,
quadratic in system size, to be able to prove the so-called one
block estimate. We do not deal with the proof of the spectral gap
bound, but we remark that with the techniques of
\cite{landimsethuramanvaradhan} one can get the desired gap
estimates for a large class of systems. Our result holds for $\b
\in (0,\frac15)$. Assuming the stronger (but harder to prove)
logarithmic-Sobolev bound we could get the result for $\b \in
(0,\frac13)$.

\section{Microscopic models}
\label{section:hyp_models}

We consider the family of microscopic models investigated in
\cite{tothvalko2}. We go over the definitions and the important
properties, for the details we refer the reader to the original
paper. There are several concrete examples introduced in
\cite{tothvalko2}, we do not list them here.

\subsection{State space, conserved quantities, generator}
\label{subs:hyp_statespace}

Throughout this paper we denote by $\Tn$ the discrete tori
$\Z/n\Z$, $n\in\N$, and by $\T$ the continuous torus $\R/\Z$. We
will denote the local spin state by $\Om$, we only consider the
case when $\Om$ is finite. The state space of the interacting
particle system is
\begin{equation*}
\Omn:=\Om^{\Tn}.
\end{equation*}
Configurations will be denoted
\begin{equation*}
\uo:=(\omega_j)_{j\in \Tn}\in\Omn,
\end{equation*}
%
%For sake of simplicity we consider discrete (integer valued)
%conserved quantities only.
The two conserved quantities are denoted by
\beqs \zeta: \Om\rightarrow \Z,\quad \eta: \Om\rightarrow \Z,
\eeqs
we also use the notations $\zeta_j=\zeta(\omega_j),\,
\eta_j=\eta(\omega_j)$.
%This means that the sums $\sum_j\zeta_j$
%and $\sum_j\eta_j$ are conserved by the dynamics.
We assume that the conserved quantities are different and
non-trivial, i.e. the functions $\zeta ,\eta$ and the constant
function 1 on $\Om$ are linearly independent.
%
%%%%%%%%%%%%%%%%%%%%%%%%%%%%%%%%%%%%%%%%%%%%%%%%%%%%%%%%%%%%%%%%%%%%%
%
%\subsection{Rate function, infinitesimal generator}
%\label{subs:hyp_rates}

We consider the \emph{rate function} $r: \Om\times \Om\times
\Om\times \Om\rightarrow \R_+$. The dynamics of the system
consists of elementary jumps effecting nearest neighbor spins,
$(\omega_j,\omega_{j+1})\longrightarrow(\omega'_j,\omega'_{j+1})$,
performed with rate
$r(\omega_j,\omega_{j+1};\omega'_j,\omega'_{j+1})$.

We require that the rate function  $r$ satisfy the following
conditions:

\begin{enumerate}[(A)]

\item \label{cond:hyp_cons} If
$r(\omega_1,\omega_2;\omega'_1,\omega'_2)>0$ then
\beq \label{eq:hyp_conserved}
\begin{array}{rcl}
\zeta(\omega_1)+\zeta(\omega_2) &=&
\zeta(\omega'_1)+\zeta(\omega'_2),
\\
\eta(\omega_1)+\eta(\omega_2) &=& \eta(\omega'_1)+\eta(\omega'_2).
\end{array}
\eeq This means that $\zeta$ and $\eta$ are indeed conserved
quantities.

\item \label{cond:hyp_irred} For every $Z\in [n \min \zeta,n \max
\zeta]\cap\Z, N\in
[n \min \eta,n \max \eta]\cap \Z$ the set %
\beqs \Omn_{Z,N}:= \left\{ \uo\in\Omn:\sum_{j\in\Tn}\zeta_j=Z,
\sum_{j\in\Tn}\eta_j=N \right\}\eeqs is an irreducible component
of $\Omn$, i.e. if $\uo, \uo'\in \Omn_{Z,N}$ then there exists a
series of elementary jumps with positive rates transforming $\uo$
into $\uo'$. This ensures that there are no hidden conservation
laws.

%If $ \omega_1,\omega_2,\omega'_1,\omega'_2 \in S$ are spins
%satisfying equations (\ref{eq:hyp_conserved}) then there exists a
%series of
%\\
%EZ NEM TELJESEN NYILVANVALO NEKEM
% pairs of spins
% $x_i,y_i$ ($0\leq i \leq n$) with $x_0=x,y_0=y$ and
% $x_n=x',y_n=y'$ or $x_0=x',y_0=y'$ and $x_n=x,y_n=y$ such that
% $r(x_i,y_i,x_{i+1},y_{i+1})>0$ for every $0\leq i<n$.

\item \label{cond:hyp_staci} There exists a probability measure
$\pi$ on $\Om$ which puts positive mass on each element of $\Om$
and for any $\omega_1$, $\omega_2$, $\omega'_1$, $\omega'_2$ $\in
\Om$
\beqs
Q(\omega_1,\omega_2)+Q(\omega_2,\omega_3)+Q(\omega_3,\omega_1)=0,\eeqs
where \beqs Q(\omega_1,\omega_2):=\sum_{\omega'_1,\omega'_2\in
\Om} \left\{
\frac{\pi(\omega'_1)\pi(\omega'_2)}{\pi(\omega_1)\pi(\omega_2)}
r(\omega'_1,\omega'_2;\omega_1,\omega_2)-r(\omega_1,\omega_2;\omega'_1,\omega'_2)
\right\}. \eeqs This condition will imply that the measure
$\prod_{j \in \Tn} \pi$ is stationary for our process on $\Omn$.
\setcounter{aux}{\value{enumi}}
\end{enumerate}

For a precise formulation of the infinitesimal generator on $\Omn$
we first define the map $\Theta_{j}^{\omega',\omega''} :\Omn
\rightarrow \Omn$ for every $\omega',\omega''\in \Om$, $j \in
\Tn$:
\beqs \left( \Theta_{j}^{\omega',\omega''} \uo \right)_{i} =
\left\{
\begin{array}{lcl}
\omega'\quad&\text{ if }& i=j
\\
\omega''\quad&\text{ if }& i=j+1
\\
\omega_i\quad&\text{ if }& i\not=j,j+1.
\end{array}
\right. \eeqs
The infinitesimal generator of the process defined on $\Omn$ is
\[
\Ln f(\uo)=\sum_{j\in \Tn} \sum_{\omega',\omega''\in \Om}
r(\omega_j,\omega_{j+1};\omega',\omega'') (f(
\Theta_{j}^{\omega',\omega''} \uo)-f(\uo)).
\]
We denote by $\Xn_t$ the Markov process on the state space $\Omn$
with infinitesimal generator $\Ln$.

%\medskip
%\noindent {\bf Remarks:}
%
%\begin{enumerate}[(1)]
%
%\item Condition (\ref{cond:hyp_cons}) implies that $\sum_j\zeta_j$ and
%$\sum_j\eta_j$ are indeed conserved quantities of the dynamics,
%while condition (\ref{cond:hyp_irred}) ensures that there are no other
%hidden conservation laws.
%
%\item Condition (\ref{cond:hyp_irred}) is somewhat implicit. It seems
%to us that it is far not trivial (if not impossible) to formulate
%explicit conditions involving the rate functions which would be
%necessary and sufficient for irreducibility. %However, in the
%%concrete examples treated in section \ref{section:hyp_examples} one
%%can easily check that irreducibility holds.
%
%\item Conditions (\ref{cond:hyp_cons}), (\ref{cond:hyp_irred}) and
%(\ref{cond:hyp_staci}) determine the measure $\pi(\omega)$ \emph{up to
%an exponential distortion}, that is the probability measures
%satisfying these conditions are of the form (\ref{eq:hyp_gibbs1}) of
%the next subsection.
%
%\item Condition (\ref{cond:hyp_staci}) implies that the stationary
%measures of the  process $\Xn_t$  are computable and have the
%structure required for hydrodynamic behaviour. See the next
%subsection for details. %Another consequence of these conditions is
%%Lemma \ref{lemma:hyp_OR} which turns out to be of  crucial importance
%%for hydrodynamic behaviour.
%
%\end{enumerate}

%%%%%%%%%%%%%%%%%%%%%%%%%%%%%%%%%%%%%%%%%%%%%%%%%%%%%%%%%%%%%%%%%%%%

\subsection{Stationary measures}
\label{subs:hyp_measures}

For every $\theta, \tau \in \R$ let $G(\theta,\tau)$ be the moment
generating function defined below:
\beqs G(\theta,\tau):=\log \sum_{\omega\in \Om} e^{\theta
\zeta(\omega)+\tau \eta(\omega)} \pi(\omega). \eeqs
%
%In thermodynamic terms $G(\theta,\tau)$ corresponds to the Gibbs
%free energy, see \cite{reichl}.
We define the probability measures
\beq \label{eq:hyp_gibbs1} \pi_{\theta, \tau}(\omega):=\pi(\omega)
\exp(\theta \zeta(\omega)+\tau \eta(\omega)-G(\theta,\tau)) \eeq
on $\Om$. Using condition (\ref{cond:hyp_staci}), by very similar
considerations as in \cite{balazs}, \cite{cocozza},
\cite{rezakhanlou} or \cite{tothvalko1} one can show that for any
$\theta,\tau \in \R$ the product measure
\beqs \pin_{\theta,\tau} :=\prod_{j \in \Tn} \pi_{\theta,\tau}
\eeqs
is stationary for the Markov process $\Xn_t$ on $\Omn$ with
infinitesimal generator $\Ln$. We will refer to these measures as
the \emph{canonical} measures. Since $\sum_j \zeta_j$ and $\sum_j
\eta_j$ are conserved, the  canonical measures on $\Omn$ are not
ergodic. The conditioned measures defined on $\Omn_{Z,N}$ by:
\beqs \pin_{Z,N}(\uo):= \pin_{\theta,\tau}\left(\uo \left|\sum_j
\zeta_j=Z,\, \sum_j \eta_j=N\right.\right) =
\frac{\pin_{\theta,\tau}(\uo)\ind\{\uo\in\Omn_{Z,N}\}}
{\pin_{\theta,\tau}(\Omn_{Z,N})} \eeqs
are also stationary and due to condition (\ref{cond:hyp_irred})
satisfied by the rate functions  they are also ergodic. We shall
call these measures the \emph{microcanonical measures} of our
system. It is easy to see that the measure $\pin_{Z,N}$ does not
depend on the values of $\theta, \tau$.

%The elementary movements of the reversed stationary process are
%$(\omega_{j-1},\omega_j)$ $\longrightarrow$
%$(\omega'_{j-1},\omega'_j)$ with rate $
%r(\omega_{j},\omega_{j-1};\omega'_j,\omega'_{j-1})$. The reversed
%generator  is
%%
%\beqs \Lna f(\uo)=\sum_{j\in \Tn} \sum_{\omega',\omega'' \in S}
%r(\omega_j,\omega_{j-1};\omega'',\omega')
%(f(\Theta_{j-1}^{\omega',\omega''} \uo)-f(\uo)). \eeqs
%%
%This is the adjoint of the operator $\Ln$ with respect to all
%microcanonical (and  canonical) measures. I.e. the reversed
%process is the same for any $\pin_{\theta,\tau}$ or $\pin_{Z,N}$.

\subsection{Expectations, fluxes}
\label{subs:hyp_expectations}

Expectation, variance, covariance with respect to the measures
$\pin_{\theta,\tau}$ will be denoted by
$\expect_{\theta,\tau}(.)$, $\var_{\theta, \tau}(.)$,
$\cov_{\theta, \tau}(.)$.

We compute the expectations of the conserved quantities with
respect  to the  canonical measures, as functions of the
parameters $\theta$ and $\tau$:
\beqs u(\theta,\tau) &:=& \expect_{\theta,\tau}(\zeta) =
\sum_{\omega \in \Om} \zeta(\omega) \pi_{\theta,\tau}(\omega)  =
\partial_\theta G(\theta,\tau)=G_{\theta},
\\
v(\theta,\tau) &:=& \expect_{\theta,\tau}(\eta) = \sum_{\omega \in
\Om} \eta(\omega) \pi_{\theta,\tau}(\omega) =
\partial_\tau G(\theta,\tau)=G_{\tau}. \eeqs
We will usually note partial derivatives by using the respective
subscripts, as long as it does not cause confusion. Elementary
calculations show, that the matrix-valued function
\beqs \left(
\begin{array}{cc}
   u_{\theta} &  u_{\tau}
\\
   v_{\theta} &  v_{\tau}
\\
\end{array}
\right)=\left(
\begin{array}{cc}
   G_{\theta\theta} &  G_{\theta\tau}
\\
   G_{\theta\tau}   &  G_{\tau\tau}
\\
\end{array}
\right) =: G(\theta,\tau) \eeqs
is equal to the covariance matrix $\cov_{\theta,
\tau}(\zeta,\eta)$ and %therefore it is strictly positive definit.
as a consequence the function $(\theta,\tau)
\mapsto(u(\theta,\tau),v(\theta,\tau))$ is invertible. We denote
the inverse function by $(u,v)\mapsto(\theta(u,v),\tau(u,v))$.
Denote by $(u,v)\mapsto S(u,v)$ the convex conjugate (Legendre
transform) of the strictly convex function $(\theta,\tau)\mapsto
G(\theta,\tau)$:
\beq \label{eq:hyp_thdentropy} S(u,v):= \sup_{\theta,\tau} \big(
u\theta+v\tau-G(\theta,\tau) \big), \eeq
and
\beq\label{eq:hyp_DD} \DD&:=&\{(u,v)\in\R_+\times\R:
S(u,v)<\infty\}
\\[5pt]
\notag &=& \text{co}\{(\zeta(\omega),\eta(\omega)):\om\in \Om\},
\eeq
where $\text{co}$ stands for convex hull. In probabilistic terms:
$S(u,v)$ is the rate function for joint large deviations of
$(\sum_j\zeta_j, \sum_j\eta_j)$. If $(u,v)$ is inside $\DD$ then
we have
\beqs \theta(u,v)=S_u(u,v), \qquad \tau(u,v)=S_v(u,v). \eeqs
%
%In
%thermodynamic terms: $S(u,v)$ corresponds to the equilibrium
%thermodynamic entropy, see \cite{reichl}.
%Let
%%
%\beqs \left(
%\begin{array}{cc}
%   \theta'_{u} &  \theta'_{v}
%\\
%   \tau'_{u} &     \tau'_{v}
%\\
%\end{array}
%\right)=\left(
%\begin{array}{cc}
%   S''_{uu} &  S''_{uv}
%\\
%   S''_{uv}   &  S''_{vv}
%\\
%\end{array}
%\right) =: S''(u,v). \eeqs
%%
%It is obvious that the matrices $G''(\theta,\tau)$ and $S''(u,v)$
%are strictly positive definit and are inverse of each other:
%%
%\beq \label{eq:hyp_inverse} G''(\theta,\tau)S''(u,v)=I, \eeq
%%
%where either $(\theta,\tau)=(u(\theta,\tau),v(\theta,\tau)) $ or
%$(u,v)=(\theta(u,v),\tau(u,v))$.
With slight abuse of notation we shall denote:
$\pi_{\theta(u,v),\tau(u,v)}=:\pi_{u,v}$,
$\pin_{\theta(u,v),\tau(u,v)}=:\pin_{u,v}$,
$\expect_{\theta(u,v),\tau(u,v)}=:\expect_{u,v}$, etc. Clearly,
$\pi_{u,v}$ can be defined naturally on the boundary of $\DD$, in
that case $\pi_{u,v}$ does not puts zero weight on some of the
elements of $\Om$.

We introduce the flux of the conserved quantities. The
infinitesimal generator $\Ln$ acts on the conserved quantities as
follows:
\beqs \Ln\zeta_i= - \phi(\omega_{i},\omega_{i+1}) +
\phi(\omega_{i-1},\omega_{i}) =: -\phi_{i} + \phi_{i-1},
\\
\Ln\eta_i= - \psi(\omega_{i},\omega_{i+1}) +
\psi(\omega_{i-1},\omega_{i}) =: -\psi_{i} + \psi_{i-1}, \eeqs
where
\beq \label{eq:hyp_phipsidef}
\begin{array}{rrl}
\phi(\omega_1,\omega_2) &:=& \dsum_{\omega'_1,\omega'_2\in \Om}
r(\omega_1,\omega_2;\omega'_1,\omega'_2)
(\zeta(\omega'_2)-\zeta(\omega_2))+C_1
\\[8pt]
%&=& \dsum_{\omega'_1,\omega'_2\in \Om}
%r(\omega_1,\omega_2;\omega'_1,\omega'_2)
%(\zeta(\omega_1)-\zeta(\omega'_1))+C,
%\\[18pt]
\psi(\omega_1,\omega_2) &:=& \dsum_{\omega'_1,\omega'_2\in \Om}
r(\omega_1,\omega_2;\omega'_1,\omega'_2)
(\eta(\omega'_2)-\eta(\omega_2))+C_2
%\\[8pt]
%&=& \dsum_{\omega'_1,\omega'_2\in \Om}
%r(\omega_1,\omega_2;\omega'_1,\omega'_2)
%(\eta(\omega_1)-\eta(\omega'_1))+C.
\end{array}
\eeq
(The constants $C_1,C_2$ may be chosen arbitrarily, we will fix
them later.) We shall denote the expectations of these functions
with respect to the canonical measure $\pi^{^{_{2}}}_{u,v}$  by
\beq \label{eq:hyp_PhiPsidef} \Phi(u,v) := \expect_{u,v}
(\phi),\qquad\qquad \Psi(u,v) := \expect_{u,v} (\psi). \eeq
The following lemma was proved in \cite{tothvalko2}.
\begin{lemma}\label{lemma:hyp_onsager}
Suppose we have a particle system with two conserved quantities
and rates satisfying conditions (\ref{cond:hyp_cons}) and
(\ref{cond:hyp_staci}). Then
\[
\partial_{\theta} \Psi\left(u(\theta,\tau),v(\theta,\tau)\right)
=\partial_{\tau} \Phi\left(u(\theta,\tau),v(\theta,\tau)\right).
\]
\end{lemma}

\noindent The first derivative matrix of the fluxes $\Phi$ and
$\Psi$ (with resp.~to $u,v$) will be denoted by
\beq \label{eq:hyp_fluxderiv} D=D(u,v):= \left(
\begin{array}{cc}
\Phi_u&\Phi_v
\\
\Psi_u&\Psi_v
\end{array}
\right). \eeq
From Lemma \ref{lemma:hyp_onsager} it follows that  $D(u,v)$ is
(weakly) hyperbolic, it can be diagonalized in a real sense (see
\cite{tothvalko2}). We denote the two eigenvalues of $D$ by
$\lambda$ and $\mu$, and the corresponding right and left
eigenvectors by
$\mathbf{r}=(r_1,r_2)^{\dag}, \mb{s}=(s_1,s_2)^{\dag}$ and $\mb{l}=(l_1,l_2), \mb{m}=(m_1,m_2)$:%
\beqs%
&&D\mb{r}=\lambda \mb{r}, \quad \mb{l}D=\lambda \mb{l},\\
&&D\mb{s}=\mu \mb{s}, \quad \mb{m}D=\mu \mb{m}.
\eeqs%
Although we do not denote it explicitly, all of these are functions of $(u,v)$. We can assume%
\beqs%
\abs{\mb{r}}=\abs{\mb{s}}=1,\quad \mb{l}\cdot \mb{r}=1, \quad
\mb{m}\cdot \mb{s}=1.\eeqs The second derivatives of the
macroscopic fluxes are denoted by $\Phi'',\Psi''$, these are
symmetric two-by-two matrices depending on $(u,v)$.

\subsection{The spectral gap condition}
\label{subs:hyp_gap}

Let $l$ be a positive integer and  $(Z,N)$ integers with
$Z\in[l\min\zeta, l\max\zeta]$, $N\in[l\min\eta, l\max\eta]$.
Expectation with respect to the measure $\pi^l_{Z,N}$ is denoted
by $\expect^l_{Z,N}\big(\cdot\big)$.  For $f:\Omega^l_{Z,N}\to\R$
let
\begin{eqnarray*}
&& \hskip-3mm L^l_{Z,N}f(\uome) := \sum_{j=1}^{l-1}
\sum_{\ome',\ome''} r(\ome_j,\ome_{j+1};\ome',\ome'')
\big(f(\Theta_{j,j+1}^{\ome',\ome''}\uome)-f(\uome)\big),
\\[5pt]
&& \hskip-3mm D^l_{Z,N}(f) := \frac12 \sum_{j=1}^{l-1}
\expect^l_{Z,N} \left( \sum_{\ome',\ome''}
r(\ome_j,\ome_{j+1};\ome',\ome'')
\big(f(\Theta_{j,j+1}^{\ome',\ome''}\uome)-f(\uome)\big)^2
\right).
\end{eqnarray*}
$L^l_{Z,N}$ is the infinitesimal generator restricted to the
hyperplane $\Omega^l_{Z,N}$, and $D^l_{Z,N}$ is the Dirichlet form
associated to $L^l_{Z,N}$ (or to its symmetric part). Note, that
$L^l_{Z,N}$ is defined with \emph{free boundary conditions}.

We will assume the following additional condition on our models:

\begin{enumerate}[(A)]
\setcounter{enumi}{\value{aux}}

\item \label{cond:hyp_gap} There exists a positive constant $W$
independent of $l,Z,N$ such that for any
$f:\Omega^l_{Z,N}\rightarrow \R$ with
$\expect^l_{Z,N} f=0$ %
\beqs \expect^l_{N,Z} f^2\le W\, l^2\, D^l_{Z,N}(f).
\eeqs%
\end{enumerate}

\noindent \textbf{Remark.} \noindent \textbf{Remark.} Presumably
(\ref{cond:hyp_gap}) is true for all (or a large class of) the
models satisfying conditions
(\ref{cond:hyp_cons})-(\ref{cond:hyp_staci}), the techniques of
\cite{landimsethuramanvaradhan} should be suitable to get the
desired gap estimates. We do not know about any published results
covering our case.

%Usually (\ref{cond:hyp_gap}) is not considered as a condition,
%rather it is proved as a lemma or theorem for the actual model or
%model-class. We chose not deal with the proof to avoid making the
%paper being too lengthy and too technical. However, we note that
%the techniques of \cite{landimsethuramanvaradhan} should be
%suitable to get the desired gap estimates for a large class of
%systems, maybe for all the  models satisfying conditions
%(\ref{cond:hyp_cons})-(\ref{cond:hyp_staci}).

% As a
%general convention,  if $\delta:\Om^m\to\R$ is a local function then
%its expectation  with respect to the canonical measure
%$\pi^{^{_{m}}}_{u,v}$ is denoted by
%%
%\beqs \Delta(u,v) := \expect_{u,v} (\delta) =
%\sum_{\omega_1,\dots,\omega_m\in \Om^m}
%\delta(\omega_1,\dots,\omega_m) \pi_{u,v}(\omega_1) \cdots
%\pi_{u,v}(\omega_m). \eeqs
%%

%
%, we use the same notations: %
%\beqs%
%&&r(\om_1,\om_2;\om_1',\om_2'), \zeta, \eta, \theta, \tau,
%\pi_{\theta, \tau}\\ &&u, v, G, \Om, \psi, \phi, \Psi, \Phi, D(u,v)
% \eeqs
%
%
%%We denote the
%%second derivative matrices of the macroscopic flux functions by%
%%\beqs%
%%A(u,v)=\nabla \Phi(u,v),\quad B(u,v)=\nabla \Phi(u,v).
%%\eeqs%
%
%The following lemma was proved in \cite{tothvalko2}.
%
%\begin{lemma}\label{lemma:hyp_onsager}
%Suppose ...
%
%Then
%\[
%\partial_{\theta} \Psi\left(u(\theta,\tau),v(\theta,\tau)\right)
%=\partial_{\tau} \Phi\left(u(\theta,\tau),v(\theta,\tau)\right).
%\]
%
%\end{lemma}

\section{Perturbation of the Eulerian hdl}
\label{section:hyp_pert}

In \cite{tothvalko2} it was proved by the application of Yau's
relative entropy method, that under Eulerian scaling the local
density profiles of the conserved quantities evolve according to
the following system of partial differential equations:
\beq \label{eq:hyp_pde_euler} \left\{\begin{array}{rcl}
  \pt u+\px \Phi(u,v)&=&0
\\[5pt]
  \pt v+\px \Psi(u,v)&=&0.
\\
\end{array}
\right. \eeq%

\noindent This pde is usually a strictly hyperbolic conservation
law (i.e. $D(u,v)$ has two distinct real eigenvalues), weak
hyperbolicity follows from Lemma \ref{lemma:hyp_onsager} (see
\cite{tothvalko2}). Since the relative entropy method needs
smoothness conditions for the solution of the limiting equation,
the previous result holds only up to a finite time, till the
appearance of the first shock. We also note, that
\cite{ollavaradhanyau} gives the first major result about the
Eulerian hydrodynamic limit of multi-component hyperbolic systems,
also wit the application of Yau's method.

\subsection{Formal perturbation}
\label{subs:hyp_formal_pert}

We will investigate the hydrodynamic behavior of small
perturbations of an equilibrium point. For that we need to
understand the asymptotics of small perturbations of a constant
solution of (\ref{eq:hyp_pde_euler}). One of the perturbation
techniques is the so-called method of weakly nonlinear geometric
optics (see e.g.~\cite{dipernamajda,hunterkeller}) which gives the
following \emph{formal} result.

Fix a point $(u_0, v_0)$ in $\DD$ and suppose that this point is
strictly hyperbolic, i.e. %
\beq%
\lambda \neq \mu,
%\textup{ and}\qquad \nabla \lambda \cdot
%\mb{r}\neq 0, \quad \nabla \mu \cdot \mb{s} \neq 0
\label{eq:hyp_cond_hyp}
\eeq%
at $(u_0,v_0)$. Suppose $(u_{\veps}(t,x), v_{\veps}(t,x))$ is the
solution of the pde
(\ref{eq:hyp_pde_euler}) with initial conditions%
\beqs%
u_{\veps}(0,x)&=&u_0+\veps u^{*}(x),\\
v_{\veps}(0,x)&=&v_0+\veps v^{*}(x),
\eeqs%
where $u^{*}(x),v^{*}(x)$ are fixed $\T\mapsto \R$ smooth
functions.
Denote%
\beq%
\begin{array}{rll}
&\sigma_0(x):=\mb{l}\cdot
(u^*(x),v^*(x))^{\dag},\quad&c_{\sigma}:=\int_{\T}\sigma_0(y)dy,\\[7pt]
&\delta_0(x):=\mb{m}\cdot
(u^*(x),v^*(x))^{\dag},\quad&c_{\delta}:=\int_{\T}\delta_0(y)dy,
\end{array}\label{eq:hyp_def_sigdel}
\eeq%
and %
\beq%
\begin{array}{rll}
a_1&:=\mb{l}\cdot (\mb{r}^{\dag} \,\Phi'' \,\mb{r}, \mb{r}^{\dag}
\,\Psi''\, \mb{r})^{\dag},\quad& a_2:=\mb{l}\cdot (\mb{r}^{\dag}
\,\Phi'' \,\mb{s},
\mb{r}^{\dag} \,\Psi''\, \mb{s})^{\dag},\\[7pt]
b_1&:=\mb{m}\cdot (\mb{s}^{\dag} \,\Phi''\, \mb{s}, \mb{s}^{\dag}
\,\Psi'' \mb{s})^{\dag},\quad& b_2:=\mb{m}\cdot (\mb{r}^{\dag}
\,\Phi''\, \mb{s}, \mb{r}^{\dag} \,\Psi'' \mb{s})^{\dag},
\end{array}\label{eq:hyp_def_ab}
\eeq%
where $\mb{l},\mb{m},\mb{r},\mb{s}$ and $\Phi'',\Psi''$ are the
respective vector- and matrix-valued functions taken at
$(u_0,v_0)$.

Then, according to the formal computations of the geometric optics
method,
\beq%
\left(%
\begin{array}{c}
  u_{\veps}(t,x) \\
  v_{\veps}(t,x) \\
\end{array}%
\right)=
\left(%
\begin{array}{c}
 u_0 \\
 v_0 \\
\end{array}%
\right)+\veps \sigma(\veps t,x-\lambda t)\, \left(\begin{array}{c}
 r_1 \\
 r_2 \\
\end{array}\right)%
+\veps \delta(\veps t,x-\mu t)\, \left(\begin{array}{c}
 s_1 \\
 s_2 \\
\end{array}\right)+\Ordo(\veps^2)%
, \label{eq:hyp_approx}
\eeq%.
as $\veps \rightarrow 0$, where $\sig$ and $\delta$ are the
solutions of the following Cauchy problems:
\beq%
\left\{\begin{array}{rrll} &\pt \sigma(t,x)+\px \left(a_1
\cdot\frac{1}{2}\sigma(t,x)^2+
c_{\delta} \, a_2\, \sigma(t,x) \right)&=&0,\\[5pt]
&\sigma(0,x)&=&\sigma_0(x),
\end{array}\right.\label{eq:hyp_def_sig_t}
\eeq%
and
\beq%
\left\{\begin{array}{rrll}&\pt \delta(t,x)+\px \left(b_1
\cdot\frac{1}{2}\delta(t,x)^2+ c_{\sigma} \, b_2\, \delta(t,x)
\right) &=&0,\\[5pt]&\delta(0,x)&=&\delta_0(x).
\end{array}\right.\label{eq:hyp_def_del_t}
\eeq%

\noindent \textbf{Remarks}

\noindent\textup{1.} This result means that a small perturbation
of a constant solution of (\ref{eq:hyp_pde_euler}) is governed by
the solutions of two decoupled equations (at least, by formal
computations). If $a_1$ and $b_1$ are nonzero, then these
equations are linear transforms of the Burgers equation. Otherwise
the respective equations become linear transport equations. It is
easy to check, that $a_1\neq 0$, $b_1 \neq 0$ hold exactly when
the point $(u_0,v_0)$ is \emph{genuinely nonlinear}, i.e.
\[
\nabla \lambda \cdot \mb{r}\neq 0, \quad \nabla \mu \cdot \mb{s}
\neq 0
\]
at $(u_0,v_0)$.

\noindent {2.} The geometric optics method is based on series
expansion, thus it needs smoothness as a condition which could
only be true up to a finite time in our case. Surprisingly, this
formal method gives good approximation of the solutions even after
the shocks. In \cite{dipernamajda} the authors prove that the
equation (\ref{eq:hyp_approx}) is valid, in the sense that for any
$t>0$ the $\textup{L}_1$-norm of the difference of the two sides
is bounded by $C\, t\, \veps^2$. In fact, this result is valid for
the case if we consider the pde (\ref{eq:hyp_pde_euler}) on $\T$
(as we do), on $\R$ they have even stronger bounds.

\subsection{The main result}
\label{subs:hyp_main}

Our main theorem is a similar result on the microscopic level. We
will apply Yau's method, thus our results will hold in the regime
of smooth solutions, only up to a finite time before the first
appearance of shocks.

Suppose, that $(u_0,v_0)$ is a point in the physical domain which
is strictly hyperbolic, see (\ref{eq:hyp_cond_hyp}). Let
$u^*(x),v^*(x)$ be smooth real functions on $\T$. Define
$\sigma(t,x), \delta(t,x)$ according to (\ref{eq:hyp_def_sigdel}),
(\ref{eq:hyp_def_ab}), (\ref{eq:hyp_def_sig_t}) and
(\ref{eq:hyp_def_del_t}), and suppose that they are smooth in
$\T\times [0,T]$. Fix a small positive parameter $\beta$, and
suppose that a particle system on $\Omn$ satisfying conditions
(\ref{cond:hyp_cons})-(\ref{cond:hyp_gap}) has initial
distribution for which the density profiles of the two conserved
quantities are 'close' to the functions $u_0+n^{-\b} u^*(\cdot)$,
$v_0+n^{-\b} v^*(\cdot)$. I.e.~the profiles are a small
perturbation of the constant $(u_0,v_0)$ profile. We also assume,
that $(u_0+n^{-\b} u^*(x)$, $v_0+n^{-\b} v^*(x))\in \cal{D}$ holds
for every $x\in \T$, at least for $n>n_0$. Then, uniformly for
$0\le t \le T$, at time $n^{1+\b} t$ the respective density
profiles will be 'close' to the functions
$u_0+n^{-\b} u^{(n)}(t,\cdot), v_0+n^{-\b} v^{(n)}(t,\cdot)$, where %
\beq%
\label{eq:hyp_def_unvn}
\left(%
\begin{array}{c}
  u^{(n)}(t,x) \\
  v^{(n)}(t,x) \\
\end{array}%
\right) :=\sigma(t,x-\lambda n^{\b} t)\,\left(%
\begin{array}{c}
  r_1 \\
  r_2 \\
\end{array}%
\right)+ \delta(t,x-\mu n^{\b} t)\,
\left(%
\begin{array}{c}
  s_1 \\
  s_2 \\
\end{array}%
\right). \eeq%
%(Where $r=(r_1,r_2)^T$, $s=(s_1,s_2)^T$.)
For the precise formulation of the result we need to introduce
some additional notations.  We will denote by $\mun_t$ the true
distribution of the system at microscopic time $n^{1+\b} t$:
\beq%
\mun_t:=\mun_0 \exp\{n^{1+\b}t L^n\}\label{eq:hyp_def_munt}
\eeq%
We define the
time-dependent reference measure $\nun_t$ as%
\beq%
\nun_t:=\prod_{j \in \Tn} \pi_{u_0+n^{-\b}
u^{(n)}(t,\frac{j}{n}),v_0+n^{-\b}v^{(n)}(t,\frac{j}{n})}\label{eq:hyp_def_nunt},
\eeq%
with $u^{(n)}, v^{(n)}$ defined in (\ref{eq:hyp_def_unvn}). This
measure mimics on a microscopic level the macroscopic profiles
$u_0+n^{-\b} u^{(n)}(t,\cdot),v_0+n^{-\b} v^{(n)}(t,\cdot)$. We
also choose an absolute
reference measure %
\beq%
\pin:=\prod_{j \in \Tn} \pi_{u_0^n,v_0^n}\label{eq:hyp_def_pin},
\eeq%
which is a stationary measure of our Markov process on $\Omn$. The
point $(u_0^n,v_0^n)$ is chosen in a way that it lies
\emph{inside} the domain $\cal{D}$ and%
\beq%
\label{eq:hyp_ref_meas_cond}
\abs{u_0-u_0^n}+\abs{v_0-v_0^n}<n^{-\b}.
\eeq%
If $(u_0,v_0)$ is inside $\cal{D}$, then we may choose
$(u_0^n,v_0^n)=(u_0,v_0)$. By choosing $(u_0^n,v_0^n)$ inside
$\DD$ we get that any probability measure on $\Omn$ is absolutely
continuous with respect to $\pin$. Condition
(\ref{eq:hyp_ref_meas_cond}) ensures that $\pin$ is 'close enough'
to $\mun_t$ in entropy sense, uniformly in $t$.

\begin{theorem*}\label{thm:hyp_main}
Let $\b\in(0,\frac15)$ be fixed. Under the stated conditions, if
\beq%
H(\mun_0|\nun_0)=o(n^{1-2\b}),\label{eq:hyp_entr_0}
\eeq%
then%
\beq%
H(\mun_t|\nun_t)=o(n^{1-2\b}), \label{eq:hyp_entr_t}\eeq%
uniformly for $0\leq t\leq T$.
\end{theorem*}

\noindent The following corollary is a simple consequence of the
Theorem and the entropy inequality.

\begin{corollary*}\label{cor:hyp_main}
Assume the conditions of Theorem \ref{thm:hyp_main}. Let $g:\T\to\R$ be a test
function. Then for any $t\in[0,T]$%
\beqs%
&&\hskip-5mm\abs{n^{-1+\b} \sum_{j\in \Tn} g(\frac{j}{n})
\left(\zeta_j(n^{1+\b}t)-u_0\right)- \int_{\T}
g(x)\left(\sigma(t,x-\lambda n^{\b} t)
  r_1 + \delta(t,x-\mu n^{\b} t)
  s_1 \right)dx }\buildrel\prob\over\rightarrow 0,\\
  &&\hskip-5mm\abs{n^{-1+\b} \sum_{j\in \Tn} g(\frac{j}{n})
\left(\eta_j(n^{1+\b}t)-v_0\right)- \int_{\T}
g(x)\left(\sigma(t,x-\lambda n^{\b} t)
  r_2 + \delta(t,x-\mu n^{\b} t)
  s_2 \right)dx }\buildrel\prob\over\rightarrow 0.
\eeqs%

\end{corollary*}

\noindent\textbf{Remarks.}

\noindent{1.} The Theorem states that if the initial distribution
of the system is 'close' to $\nun_0$ in relative entropy sense
then at time $n^{1+\b}t$ it will be close to $\nun_t$. The fact,
that 'close' should mean $o(n^{1-2\b})$ can be easily explained,
see e.g.~\cite{tothvalko1} or \cite{tothvalko2}.

\noindent{2.} If instead of condition (\ref{cond:hyp_gap}) we
assume a similar uniform bound on the logarithmic-Sobolev constant
then our Theorem is valid for $\b\in (0,\frac13)$.

%\subsection{Main result}

\section{Proof}

%\subsection{Preparatory...}
We will assume, that%
\beq (u_0,v_0)=(0,0),\qquad
\Phi_v(0,0)=\Psi_u(0,0)=0.\label{eq:hyp_simpl}
\eeq%
It is easy to see, that we can always reduce the general case to
get (\ref{eq:hyp_simpl}), via some suitable linear transformations
on $(\zeta, \eta)$. Also, with the proper choice of the constants
in the definition (\ref{eq:hyp_phipsidef}) we can set
\beq%
\Phi(0,0)=\Psi(0,0)=0.\label{eq:hyp_simpl_}
\eeq%
Assumptions (\ref{eq:hyp_simpl}) imply,
that%
\beq%
D=\left(%
\begin{array}{cc}
  \lambda & 0 \\
  0 & \mu \\
\end{array}%
\right),\label{eq:hyp_D0} \qquad \mb{l}=\mb{r}^{\dag}=(1,0),\qquad
\mb{m}=\mb{s}^{\dag}=(0,1),
\eeq%
%kell???
%
%Also, via Lemma \ref{lemma:hyp_onsager}
%\beq%
%&&\tau_u(0,0)=\theta_v(0,0)=0,\notag\\&&u_{\tau}(0,0)=v_{\theta}(0,0)=0.
%\eeq
%
%
and %
\beq%
 u^{(n)}(t,x)= \sigma(t,x-\lambda n^{\b} t),\qquad
\label{eq:hyp_unvn_0}%
v^{(n)}(t,x)= \delta(t,x-\mu n^{\b} t).
\eeq%
We introduce the notations%
\beq%
\Phi''=\left(%
\begin{array}{cc}
  \Phi_{uu} & \Phi_{uv} \\
  \Phi_{vu} & \Phi_{vv} \\
\end{array}%
\right)=:
\left(%
\begin{array}{cc}
  a_1 & a_2 \\
  a_2 & a_3 \\
\end{array}%
\right),\qquad
\Psi''=\left(%
\begin{array}{cc}
  \Psi_{uu} & \Psi_{uv} \\
  \Psi_{vu} & \Psi_{vv} \\
\end{array}%
\right)=:
\left(%
\begin{array}{cc}
  b_3 & b_2 \\
  b_2 & b_1 \\
\end{array}%
\right).\label{eq:hyp_ab_0}
\eeq%
Clearly, these definitions agree with the definition
(\ref{eq:hyp_def_ab}).

%Then, from (\ref{eq:hyp_def_ab}) and (\ref{eq:hyp_ab_0}) we have
%\[
%a=a_1,\qquad b=b_3.
%\]
We define the functions $\ols(t,x_1,x_2), \old(t,x_1,x_2)$ as%
\beq%
%\begin{array}{rlc}
\ols(t,x_1,x_2)&:=&\frac{1}{\lambda-\mu}\left(a_2 \, \sigma(t,x_1)
\delta(t,x_2)+a_2 \, \sigma_x (t,x_1) \int_0^{x_2}
(\delta(t,z)-c_{\delta})dz +\frac{a_3}{2}
\,\delta(t,x_2)^2\right)\notag
\\[-3pt]
&&\label{eq:hyp_def_olsold}
\\
\old(t,x_1,x_2)&:=&\frac{1}{\mu-\lambda}\left(b_2 \, \sigma(t,x_1)
\delta(t,x_2)+b_2 \, \delta_x (t,x_2) \int_0^{x_1}
(\sigma(t,z)-c_{\sigma})dz +\frac{b_3}{2}
\,\sigma(t,x_1)^2\right)\notag
%\end{array}
\eeq%

\noindent The defining partial differential equations
(\ref{eq:hyp_def_sig_t}), (\ref{eq:hyp_def_del_t}) of the
functions $\sigma, \delta$ are conservation laws, thus for any
$0\le t \le T$:
\[
\int_{\T}\sigma(z,t)dz=c_{\sigma},\quad
\int_{\T}\delta(z,t)dz=c_{\delta}.
\]
From that it follows that $\ols, \old$ are well-defined smooth
functions on $[0,T]\times \T \times \T$ (i.e.~periodic in $x_1$
and $x_2$) with bounded derivatives.

%?????????????Kell ide?
%
%In \cite{tothvalko2} it was proved using ... that the following
%'Onsager'-relation holds: %
%\beq%
%\partial_{\theta}
%\Psi(u(\theta,\tau),v(\theta,\tau))=\partial_{\tau}
%\Phi(u(\theta,\tau),v(\theta,\tau))\label{eq:hyp_Onsager}.
%\eeq%
%Applying this at $(\theta,\tau)=(0,0)$ and using (\ref{eq:hyp_D0}) we readily get%
%\beq
%&&\tau_u(0,0)=\theta_v(0,0)=0,\notag\\&&u_{\tau}(0,0)=v_{\theta}(0,0)=0.
%\eeq

\subsection{Changing the time-dependent reference measure}
The usual way to prove a result like Theorem \ref{thm:hyp_main} is
to
get a Grönwall-type estimate on $H(\mun_t|\nun_t)$:%
\beqs%
H(\mun_t|\nun_t)-H(\mun_0|\nun_0)\leq C \int_0^t
H(\mun_s|\nun_s)+o(n^{1-2\b}),
\eeqs%
via bounding the derivative $\pt H(\mun_t|\nun_t)$. We will use a
slightly different approach, by proving a similar estimate for
$H(\mun_t|\nuun_t)$:%
\beq\label{eq:hyp_gronwall}%%
H(\mun_t|\nuun_t)-H(\mun_0|\nuun_0)\leq C \int_0^t
H(\mun_s|\nuun_s)+o(n^{1-2\b}).
\eeq%
Here
\beq%
\nuun_t:=\prod_{j \in \Tn}
\pi_{n^{-\b}\wt{u}^{(n)}(t,\frac{j}{n}),n^{-\b}
\wt{v}^{(n)}(t,\frac{j}{n})}\label{eq:hyp_def_nuunt}
\eeq%
and $\wt{u}^{(n)},\wt{v}^{(n)}$ are smooth functions defined
as%
\beq%
\wt{u}^{(n)}(x,t)&:=&u^{(n)}(x,t)+n^{-\b}\, \ols(t,x-\lambda
n^{\b}t,x-\mu
n^{\b}t)\notag\\[5pt]&=& \sigma(t,x-\lambda n^{\b}
t)+n^{-\b}\, \ols(t,x-\lambda n^{\b}t,x-\mu n^{\b}t)
,\notag \\[5pt]
\wt{v}^{(n)}(x,t)&:=&v^{(n)}(x,t)+n^{-\b}\, \old(t,x-\lambda
n^{\b}t,x-\mu
n^{\b}t),\\[5pt]
&=& \delta(t,x-\mu n^{\b} t)+n^{-\b}\, \old(t,x-\lambda
n^{\b}t,x-\mu n^{\b}t) .\notag
\eeq%
Because of Lemma \ref{lemma:hyp_switch} and condition
(\ref{eq:hyp_entr_0}) we have $H(\mun_0|\nuun_0)=o(n^{1-2\b})$ and
therefore from (\ref{eq:hyp_gronwall})
\[H(\mun_t|\nuun_t)=o(n^{1-2\b})\]
will follow uniformly for $0\leq t \leq T$. Using Lemma
\ref{lemma:hyp_switch} again we get Theorem \ref{thm:hyp_main}.

%%%%%%%%%%%%%%%%%%%%%%%%%%%%%%%%%%%%%%%%%%%%
\begin{lemma}\label{lemma:hyp_switch}
Let $\mun_t,\nun_t,\nuun_t$ be the measures defined as before,
with $t\in[0,T]$. Then
\[
H(\mun_t|\nun_t)=o(n^{1-2\b}) \Longleftrightarrow
H(\mun_t|\nuun_t)=o(n^{1-2\b}).
\]
\end{lemma}

\begin{proof}
We start with
\beq%
 H(\mun_t|\nun_t)-H(\mun_t|\nuun_t)=-\int_{\Omn}
\log \frac{d\nun_t}{d\nuun_t} d\mun_t.\label{eq:hyp_alap}
\eeq%
By subsections \ref{subs:hyp_measures} and
\ref{subs:hyp_expectations} we
can calculate that%
\beqs%
\log \frac{d\nun_t}{d\nuun_t}(\uo)&=&\sum_{j\in \Tn } \left\{
\phantom{+}\left(\theta(n^{-\b} u^{(n)},n^{-\b} v^{(n)})
-\theta(n^{-\b} \uun,n^{-\b}
\vvn)\right)\zeta_j\right.\\&&\hskip10mm+\left(\tau(n^{-\b}
u^{(n)},n^{-\b} v^{(n)}) -\tau(n^{-\b} \uun,n^{-\b}
\vvn)\right)\eta_j\\[5pt]&&\hskip10mm-G\left(\theta(n^{-\b} u^{(n)},n^{-\b}
v^{(n)}),\tau(n^{-\b} u^{(n)},n^{-\b}
v^{(n)})\right)\\[5pt]&&\left.\hskip10mm+G\left(\theta(n^{-\b} \uun,n^{-\b}
\vvn),\tau(n^{-\b} \uun,n^{-\b} \vvn)\right) \right\},
\eeqs%
where, for typographical reasons, we omitted the arguments
$(t,\frac{j}{n})$ from the
functions $u^{(n)},v^{(n)},\uun,\vvn$.\\
From the previous expression via power-series expansion:%
\beqs%
\log \frac{d\nun_t}{d\nuun_t}(\uo)&\le& \Ordo(n^{1-3\b})+C
n^{-2\b} \sum_{j\in \Tn}
\left(\abs{\zeta_j-u^{(n)}(t,\frac{j}{n})}+\abs{\eta_j-v^{(n)}(t,\frac{j}{n})}\right),\\
&=&\Ordo(n^{1-3\b})+C n^{-2\b} \sum_{j\in \Tn}
\left(\abs{\zeta_j-\uun(t,\frac{j}{n})}+\abs{\eta_j-\vvn(t,\frac{j}{n})}\right),
\eeqs%
with uniform error terms. Using this with (\ref{eq:hyp_alap}) and
the entropy inequality with respect to $\nun_t$ and $\nuun_t$ the
lemma follows.
\end{proof}

We also note, that applying the same arguments as in the proof of
Lemma \ref{lemma:hyp_switch}  we get that from the condition
(\ref{eq:hyp_entr_0}) \[H(\mun_0|\pin)=\Ordo(n^{1-2\b})\] follows.
Since $\pin$ is a stationary measure,%
\beq H(\mun_t|\pin)\leq  H(\mun_0|\pin)=\Ordo(n^{1-2\b})
\label{eq:hyp_entr_bound_abs}\eeq for all $t\geq 0$.

The proof of the following lemma is a simple application of the
entropy inequality with the entropy bound
(\ref{eq:hyp_entr_bound_abs}). Mind that because of
(\ref{eq:hyp_simpl}) and (\ref{eq:hyp_simpl_})%
\[
\expect_{\pin}\zeta=\expect_{\pin}\eta=\expect_{\pin}\phi=\expect_{\pin}\psi=0.
\]

\begin{lemma}\label{lemma:hyp_entropy}
Suppose $b_1,b_2,\dots$ are real numbers with $\abs{b_j}\le1$ and
$\xi_j$ stands for either of $\eta_j$, $\zeta_j$, $\psi_j$ or
$\phi_j$. Then
\beqs%
&&\int_{\Omn} \frac{1}{n}\sum_{j\in\Tn} b_j \xi_j \,d\mun_t\le C
n^{-\b}. \eeqs%
with an absolute constant $C$.
\end{lemma}

In the rest of the paper we prove inequality
(\ref{eq:hyp_gronwall}).

%%%%%%%%%%%%%%%%%%%%%%%%%%%%%%%%%%%%%%%%%%%%%%%%%%%%%%%%%%%%%%%%

\subsection{Preparatory computations}

We define
\beq%
\tetn(t,x)&:=&n^{\b} \theta(n^{-\b} \uun(t,x),n^{-\b} \vvn(t,x)),\notag\\
\taun(t,x)&:=&n^{\b} \tau(n^{-\b} \uun(t,x),n^{-\b} \vvn(t,x)),
\label{eq:hyp_def_tettaun}\\
\theta^{n}_0&:=&n^{\b} \theta(n^{-\b} u^n_0,n^{-\b} v_0^n),\notag\\
\tau^{n}_0&:=&n^{\b} \tau(n^{-\b} u^n_0,n^{-\b} v_0^n).\notag
\eeq%
It is easy to check, that the partial derivatives $\partial_x
\tetn(t,x),\partial_x \taun(t,x),\partial_x^2
\tetn(t,x),\partial_x^2 \taun(t,x)$ are uniformly bounded in
$[0,T]\times \T$. From subsection \ref{subs:hyp_measures} we have
\beq%
\ffnt&:=&\frac{d\nuunt}{d\pin}\notag\\&=&\exp \sum_{j\in \Tn }
\left\{ n^{-\b} (\tetn(t,\frac{j}{n})-\theta^n_0) \zeta_j+n^{-\b}
(\taun(t,\frac{j}{n})-\tau^n_0)
\eta_j\right.\label{eq:hyp_fnt}\\&&\hskip15mm\left.-G\left(n^{-\b}
\tetn(t,\frac{j}{n}),n^{-\b}
\taun(t,\frac{j}{n})\right)+G\left(n^{-\b} \tetn_0,n^{-\b} \taun_0
\right) \right\}.\notag
\eeq%
Differentiating the identity
\[
H(\munt|\nuun_t)-H(\munt|\pin)=-\int_{\Omn} \log{\ffnt} \munt
\]
and noting that $\pt H(\munt|\pin)\leq 0$ we get the following
bound on $\pt H(\munt|\nuun_t)$:%
\beqs%
n^{2\b-1}\pt H(\munt|\nuun_t)\leq -\int_{\Omn} \left(n^{3\b}
\Ln\log{\ffnt}+n^{-1+2\b} \pt \log{\ffnt}\right) d\munt.
\eeqs%
Integrating with respect to the time:
\beq%
n^{2\b-1}(H(\munt|\nuun_t)-H(\mun_0|\nuun_0))\leq
-\int_0^t\int_{\Omn} \left(n^{3\b} \Ln\log{\ffn_s}+n^{-1+2\b} \pt
\log{\ffn_s}\right) d\mun_s dt\label{eq:hyp_initial}.
\eeq%
We estimate the two terms on the right-hand side separately in the
next two subsections.

\subsection{Estimating the first term of (\ref{eq:hyp_initial})}

From the definitions
\beq%
\hskip-15mm\notag n^{3\b}
\Ln\log{\ffnt}(\uo)&=&n^{-1+2\b}\sum_{j\in \Tn} \left(
\phi_j-\Phi(n^{-\b} \uun_j,n^{-\b}
\vvn_j)\right)\px \tetn(t,\frac{j}{n})\\
\label{eq:hyp_Lf}&&\phantom{MMMM}+\left( \psi_j-\Psi(n^{-\b}
\uun_j,n^{-\b} \vvn_j)\right)\px \taun(t,\frac{j}{n})
\\&&+\Err_1^n(t,\uo)+\Err_2^n(t),\notag
\eeq%
where%
\beq%
\Err_1^n(t,\uo)&:=&n^{-1+2\b}\sum_{j\in \Tn}\left\{  \phi_j
\left(\gradn \tetn(t,\frac{j}{n}) -\px
\tetn(t,\frac{j}{n})\right)\right. \notag
\\&&
\hskip20mm+\left.\psi_j\left(\gradn
\taun(t,\frac{j}{n}) -\px \taun(t,\frac{j}{n})\right)\right\},\label{eq:hyp_err_1}\\
\Err_2^n(t)&:=&n^{-1+2\b}\sum_{j\in \Tn} \left\{\Phi\left(n^{-\b}
\uun_j,n^{-\b} \vvn_j\right)\px \tetn(t,\frac{j}{n})\right.\notag
\\&&\hskip20mm+\left.\Psi\left(n^{-\b} \uun_j,n^{-\b} \vvn_j\right)\px
\taun(t,\frac{j}{n})\right\}.\label{eq:hyp_err_2}
\eeq%
We used the (slightly abused) shorthanded notations
\beqs%
\uun_j=\uun(t,\frac{j}{n}),\hskip20mm \vvn_j=\vvn(t,\frac{j}{n}),
\eeqs%
and $\gradn$ denotes the discrete gradient:
\[
\gradn g(x):=n\left(g(x+\frac1n)-g(x)\right).
\]
Using the smoothness of $\tetn$ and $\taun$ and Lemma
\ref{lemma:hyp_entropy} the expectation of the first error term
can be easily estimated:
\beq%
\int_{\Omn}\abs{\Err_1^n(t,\uo)}d\mun_t=
\Ordo(n^{-1+\b}).\label{eq:hyp_bound_err_1}
\eeq%
Using Lemma \ref{lemma:hyp_onsager} it is easy to see that there
exists
a smooth function $U(u,v)$ such that%
\beqs%
\partial_{\theta} U\left(u(\theta,\tau),v(\theta,\tau)\right)
=\Phi\left(u(\theta,\tau),v(\theta,\tau)\right),\qquad
\partial_{\tau} U\left(u(\theta,\tau),v(\theta,\tau)\right)
=\Psi\left(u(\theta,\tau),v(\theta,\tau)\right). \eeqs%
Thus $\Err_2^n(t)$ takes the form:
\[
\Err_2^n(t,\uo)=n^{3\b-1} \sum_{j\in \Tn} \px
U\left(n^{-\b}\uun(t,\frac{j}{n}),n^{-\b}\vvn(t,\frac{j}{n})\right),
\]
from which
\beq%
\Err_2(t)=\Ordo(n^{-1+2\b}),\label{eq:hyp_bound_err_2}
\eeq%
uniformly for $t\in[0,T]$. From previous bounds  we have
\beq%
&&\hskip-15mm\notag \int_0^t \int_{\Omn}n^{3\b}
\Ln\log{\ffn_s}(\uo)d\mun_s ds=\\
&&\hskip15mm \notag n^{-1+2\b}\int_0^t\int_{\Omn}\sum_{j\in \Tn}
\left( \phi_j-\Phi(n^{-\b} \uun_j,n^{-\b}
\vvn_j)\right)\px \tetn(s,\frac{j}{n})\\
\label{eq:hyp_Lf1}&&\hskip43mm+\left( \psi_j-\Psi(n^{-\b}
\uun_j,n^{-\b} \vvn_j)\right)\px \taun(s,\frac{j}{n})d\mun_t ds
\\&&\hskip15mm+\Ordo(n^{-1+2\b})\notag
\eeq%
In the next step we introduce the block averages.  We will denote
the block size with $l=l(n)$, it will be large microscopically,
but small on the macroscopic scale. In the first computations we
only assume $l\gg n^{2\b}$, the exact order of $l$ will only be
determined at the end of the proof, after collecting all the error
terms. For a local function $\kappa_j$ ($j\in\Tn$) we define its
block average with
\[
\kappa^l_j:=\frac1{l}\sum_{i=0}^{l-1}\kappa_{j+i}.
\]
By partial summation for a smooth function $\rho(x):\T\mapsto \R$
we have
\[
\abs{\sum_{j\in\Tn}\kappa_j\,
\rho(\frac{j}{n})-\sum_{j\in\Tn}\kappa_j^l \,\rho(\frac{j}{n})}\le
 \|\px \rho\|_{\infty} \, \abs{\sum_{j\in\Tn}\kappa_j} \frac{l}{n}
\]
Using this with Lemma \ref{lemma:hyp_entropy} we can replace
$\phi_j,\psi_j$ in (\ref{eq:hyp_Lf1}) with the the respective
block averages:
\beq%
&&\hskip-15mm\notag \int_0^t \int_{\Omn}n^{3\b}
\Ln\log{\ffn_s}(\uo)d\mun_s ds=\\
&&\hskip15mm \notag n^{-1+2\b}\int_0^t\int_{\Omn}\sum_{j\in \Tn}
\left\{\left( \phi_j^l-\Phi(n^{-\b} \uun_j,n^{-\b}
\vvn_j)\right)\px \tetn(s,\frac{j}{n})\right.\\
\label{eq:hyp_Lf2}&&\hskip45mm\left.+\left( \psi_j^l-\Psi(n^{-\b}
\uun_j,n^{-\b} \vvn_j)\right)\px
\taun(s,\frac{j}{n})\right\}d\mun_t ds
\\&&\hskip15mm+\Ordo(n^{\b-1}l).\notag
\eeq%
Finally, using Lemma \ref{lemma:hyp_obe} (the one-block estimate),
we replace the block averages $\phi^l_j, \psi_j^l$ by their 'local
equilibrium value': $\Phi(\zeta^l_j,\eta_j^l)$ and
$\Psi(\zeta^l_j,\eta_j^l)$, respectively:%
\beq%
&&\hskip-15mm\notag \int_0^t \int_{\Omn}n^{3\b}
\Ln\log{\ffn_s}(\uo)d\mun_s ds=\\
&&\hskip10mm \notag n^{-1+2\b}\int_0^t\int_{\Omn}\sum_{j\in \Tn}
\left(\Phi(\zeta^l_j,\eta_j^l)-\Phi(n^{-\b} \uun_j,n^{-\b}
\vvn_j)\right)\px \tetn(s,\frac{j}{n})\\
\label{eq:hyp_Lf3}&&\hskip38mm+\left(\Psi(\zeta^l_j,\eta_j^l)-\Psi(n^{-\b}
\uun_j,n^{-\b} \vvn_j)\right)\px \taun(s,\frac{j}{n})d\mun_t ds
\\&&\hskip10mm+\Ordo(n^{\b-1}l\vee
n^{-1-\b} l^3\vee l^{-1})\notag
\eeq%

\begin{lemma}[One block estimate]\label{lemma:hyp_obe}
\beqs%
\frac1{n}\int_0^t\int_{\Omn}\sum_{j\in
\Tn}\abs{\phi_j^l-\Phi(\zeta^l_j,\eta_j^l)}d\mun_s dt &\le& C
(n^{-1-2\b} l^3+ l^{-1}),\\\frac1{n}\int_0^t\int_{\Omn}\sum_{j\in
\Tn}\abs{\psi_j^l-\Psi(\zeta^l_j,\eta_j^l)}d\mun_s dt &\le& C
(n^{-1-3\b} l^3 + l^{-1}).
\eeqs%
\end{lemma}
The proof relies on the spectral gap condition
(\ref{cond:hyp_gap}). It uses the Feynman-Kac formula, the
Raleigh-Schrödinger perturbation technique and the 'equivalence of
ensembles' (see the Appendix of \cite{kipnislandim} for all
three). A detailed proof can be found in \cite{tothvalko1} for the
one component case which can be easily adapted for our purposes.

\noindent \textbf{Remark.} If instead of the condition
(\ref{cond:hyp_gap}) we have a similar uniform bound on the
logarithmic-Sobolev constant, then the previous lemma may be
strengthened: it holds with the bound $C (n^{-1-3\b} l^2 +
l^{-1}).$

\subsection{Estimating the second term of (\ref{eq:hyp_initial})}

Performing the time-derivation we obtain:%
\beq%
\hskip-10mm n^{-1+2\b} \pt \log{\ffnt}&=&\frac1{n} \sum_{j\in \Tn}
\left\{ \left(n^{\b} \zeta_j-\uun_j\right)\pt \tetn(t,\frac{j}{n})
%\right.
%\\ &&
%\phantom{\frac1{n} \sum_{j\in
%\Tn}}
%\left.
+\left(n^{\b} \eta_j-\vvn_j\right)\pt
\taun(t,\frac{j}{n})\right\}\label{eq:hyp_ptlogfnt0}
\eeq%
By the definitions of $\uun, \vvn$  and Taylor-expansion we
readily get that%
\beqs%
\pt \uun(t,x)&=&\sig_t(t,x-\lambda n^{\b} t)-\lambda n^{\b}
\sig_x(t,x-\lambda n^{\b} t)\\&&-\lambda \ols_{x_1}(t,x-\lambda
n^{\b} t,x-\mu n^{\b} t)-\mu \ols_{x_2}(t,x-\lambda n^{\b} t,x-\mu
n^{\b} t)\\&&+\Ordo(n^{-\b}),
\eeqs%
and
\beqs%
&&\hskip-15mmn^{2\b}\px \Phi\left(n^{-\b} \uun(x,t),n^{-\b}
\vvn(x,t)\right)\\&&\hskip20mm=\phantom{+}\lambda
n^{\b}\sig_x(t,x-\lambda n^{\b} t)\\&&\hskip23mm+\lambda
\ols_{x_1}(t,x-\lambda n^{\b} t,x-\mu n^{\b} t)+\lambda
\ols_{x_2}(t,x-\lambda n^{\b} t,x-\mu n^{\b} t)\\&&\hskip23mm+ \px
\left(\frac12 a_1 \sig(t,x-\lambda n^{\b} t)^2 +c_{\del} a_2
\sig(t,x-\lambda n^{\b} t)\right)\\&&\hskip23mm +\px\left(a_2
\sig(t,x-\lambda n^{\b} t) (\del(t,x-\mu n^{\b}
t)-c_{\del})+\frac12 a_3 \del(t,x-\mu n^{\b} t)^2\right)
\\&&\hskip23mm+\Ordo(n^{-\b}),
\eeqs%
with uniform error terms. ($\ols_{x_1}$ and $\ols_{x_2}$ are the
partial derivatives of $\ols(t,x_1,x_2)$ with respect to the
second and third variable.) Adding up these equations and checking
the definitions for $\sig, \del, \ols$ we see that all the
significant terms on the right hand side cancel to give:
\beq%
\pt \uun+\px \left(n^{2\b} \Phi(n^{-\b} \uun,n^{-\b} \vvn)
\right)=\Ordo(n^{-\b}).\label{eq:hyp_ptuun1}
\eeq%
Similarly,
\beq%
\pt \vvn+\px \left(n^{2\b} \Psi(n^{-\b} \uun,n^{-\b} \vvn)
\right)=\Ordo(n^{-\b}).\label{eq:hyp_ptuun2}
\eeq%
From (\ref{eq:hyp_ptuun1}) and (\ref{eq:hyp_ptuun2}):%
\beq%
\pt \tetn&=&\phantom{-}\theta_u \pt \uun+\theta_v \pt \vvn \notag \\
&=&-n^{2\b} \theta_u \px \left(\Phi(n^{-\b} \uun,n^{-\b} \vvn)
\right)-n^{2\b} \theta_v \px \left(\Psi(n^{-\b} \uun,n^{-\b} \vvn)
\right)+\Ordo(n^{-\b})\notag
\\
&=&-n^{\b} (\theta_u \Phi_u +\theta_v \Psi_u) \px \uun
-n^{\b}(\theta_u \Phi_v +\theta_v \Psi_v) \px
\vvn+\Ordo(n^{-\b})\notag\\
&=&-n^{\b} (\theta_u \Phi_u +\tau_u \Psi_u) \px \uun
-n^{\b}(\theta_v \Phi_u +\tau_v \Psi_u) \px
\vvn+\Ordo(n^{-\b})\notag\\
&=&-n^{\b} \Phi_u \px \tetn -n^{\b} \Psi_u \px
\taun+\Ordo(n^{-\b})\label{eq:hyp_ptthetn}
\eeq%
In the fourth line we used $\tau_u=\theta_v$ and Lemma
\ref{lemma:hyp_onsager}. To simplify notations, we omitted the
arguments $(t,x)$ from the functions $\tetn,\taun,\uun,\vvn$, and
the arguments $(n^{-\b}\uun,n^{-\b}\vvn)$ from all the partial
derivatives of $\theta, \tau,\Phi, \Psi$ with respect to $u,v$.
Similarly,%
\beq%
\pt \taun=-n^{\b} \Phi_v \px \tetn -n^{\b} \Psi_v \px
\taun+\Ordo(n^{-\b}).\label{eq:hyp_pttaun}
\eeq%
Hence from (\ref{eq:hyp_ptlogfnt0}):%
\beq%
\notag &&\hskip-15mm n^{-1+2\b} \pt \log{\ffnt}\\
\notag &&\hskip-8mm=-n^{-1+2\b} \sum_{j\in \Tn} \left\{ \left(
\zeta_j-n^{-\b}\uun_j\right)
\Phi_u\left(n^{-\b}\uun_j,n^{-\b}\vvn_j\right) \px
\tetn(t,\frac{j}{n}) \right.
\\&&\hskip18mm\label{eq:hyp_ptlogfnt1}
+\left( \zeta_j-n^{-\b}\uun_j\right)
\Psi_u\left(n^{-\b}\uun_j,n^{-\b}\vvn_j\right) \px
\taun(t,\frac{j}{n})
\\[3pt]\notag&&\hskip18mm
+\left( \eta_j-n^{-\b}\vvn_j\right)
\Phi_v\left(n^{-\b}\uun_j,n^{-\b}\vvn_j\right) \px
\tetn(t,\frac{j}{n}) \\[3pt]  &&\hskip18mm\left.+\left(
\eta_j-n^{-\b}\vvn_j\right)\Psi_v\left(n^{-\b}\uun_j,n^{-\b}\vvn_j\right)\px
\taun(t,\frac{j}{n})\right\}\notag\\&&+\Err_3(t,\uo),\notag
\eeq%
where%
\beq%
\Err_3(t,\uo)= \frac{1}{n}\sum_{j\in \Tn}\label{eq:hyp_err_3}
(\zeta_j-n^{-\b} \uun_j)\, b_j(t)+(\eta_j-n^{-\b} \vvn_j)\, c_j(t)
\eeq%
and $b_j(t)$ and $c_j(t)$ are uniformly bounded constants. Using
Lemma \ref{lemma:hyp_entropy} we get that
\beq%
\int_{\Omn}\abs{\Err_3(t,\uo)}d\munt=\Ordo(n^{-\b}).\label{eq:hyp_bound_err_3}
\eeq%
We can exchange $\zeta_j,\eta_j$ with their block-averages
$\zeta_j^l,\eta_j^l$ (as in the previous subsection) which (after
the time-integration) gives the following estimate:
\beq%
\notag &&\hskip-15mm \int_0^t \int_{\Omn} n^{-1+2\b} \pt \log{\ffn_s} d\mun_s ds\\
\notag &&\hskip-8mm=-\int_0^t \int_{\Omn}n^{-1+2\b} \sum_{j\in
\Tn} \left\{ \left( \zeta_j^l-n^{-\b}\uun_j\right)
\Phi_u\left(n^{-\b}\uun_j,n^{-\b}\vvn_j\right) \px
\tetn(s,\frac{j}{n}) \right.
\\&&\hskip33mm\label{eq:hyp_ptlogfnt2}
+\left( \zeta_j^l-n^{-\b}\uun_j\right)
\Psi_u\left(n^{-\b}\uun_j,n^{-\b}\vvn_j\right) \px
\taun(s,\frac{j}{n})
\\[3pt]\notag&&\hskip33mm
+\left( \eta_j^l-n^{-\b}\vvn_j\right)
\Phi_v\left(n^{-\b}\uun_j,n^{-\b}\vvn_j\right) \px
\tetn(s,\frac{j}{n}) \\[3pt]  &&\hskip33mm\left.+\left(
\eta_j^l-n^{-\b}\vvn_j\right)\Psi_v\left(n^{-\b}\uun_j,n^{-\b}\vvn_j\right)\px
\taun(s,\frac{j}{n})\right\}d\mun_s ds\notag\\&&+\Ordo(n^{-\b}\vee
n^{\b-1} l),\notag
\eeq%

%\noindent Summing up our estimates:
%\beq%
%&&\hskip-20mm-\int_{\Omn} \left(n^{3\b} \Ln\log{\ffnt}+n^{-1+2\b}
%\pt
%\log{\ffnt}\right) d\munt\notag\\
%\notag&&\hskip-10mm=- n^{-1+2\b}\int\limits_{\Omn}\sum_{j\in \Tn}
%\px \tetn(t,\frac{j}{n}) \left(
%\phi_j-\Phi(n^{-\b} \uun_j,n^{-\b} \vvn_j)\right.\\
%\notag&&\hskip40mm- \Phi_u(n^{-\b}\uun_j,n^{-\b}\vvn_j)\left(
%\zeta_j-n^{-\b}\uun_j\right)\\\notag&&\hskip40mm\left.-
%\Phi_v(n^{-\b}\uun_j,n^{-\b}\vvn_j)\left(
%\eta_j-n^{-\b}\vvn_j\right)
%\right)\\
%\label{eq:hyp_initial_est}&&\hskip20mm+\px \taun(t,\frac{j}{n}) \left(
%\psi_j-\Psi(n^{-\b} \uun_j,n^{-\b} \vvn_j)\right.\\
%\notag&&\hskip40mm- \Psi_u(n^{-\b}\uun_j,n^{-\b}\vvn_j)\left(
%\zeta_j-n^{-\b}\uun_j\right)\\\notag&&\hskip40mm\left.-
%\Psi_v(n^{-\b}\uun_j,n^{-\b}\vvn_j)\left(
%\eta_j-n^{-\b}\vvn_j\right) \right)d\mun_t\\ \notag &&
%+\Ordo(n^{-\b}\vee n^{1-2\b}).
%\eeq%

%\subsection{One block estimate}
%\label{subs:hyp_obe}\ref{}

\subsection{Block replacement}
\label{subs:hyp_block_repl}

For a function $\Upsilon(u,v)$ we denote %
\beqs \RR_{\Ups}(u_1,v_1;u_2,v_2):=\Ups(u_1,v_1)-\Ups(u_2,v_2)
-\Ups_u(u_2,v_2)(u_1-u_2)-\Ups_v(u_2,v_2)(v_1-v_2).%
\eeqs%
Collecting the estimates of the previous subsections we have%
\beq%
&&\hskip-10mm H(\mun_t|\nuun_t)-H(\mun_t|\nuun_0)\notag\\&&\notag
\hskip5mm \leq C\, n^{2\b-1} \int_0^t \int_{\Omn} \sum_{j\in \Tn}
\left(\abs{\RR_{\Phi}\left(\zeta^l_j,\eta^l_j;n^{-\b}
\uun(s,\frac{j}{n}),n^{-\b} \vvn(s,\frac{j}{n})\right)}\right.\\[-10pt]\label{eq:hyp_bl_rpl_0} \\
\notag
&&\hskip42mm\left.+\abs{\RR_{\Psi}\left(\zeta^l_j,\eta^l_j;n^{-\b}
\uun(s,\frac{j}{n}),n^{-\b} \vvn(s,\frac{j}{n})\right)}\right)d\mun_s ds\\
&&\hskip10mm+\Ordo(n^{-\b} \vee n^{\b-1} l \vee n^{-1-\b}
l^3)\notag
\eeq%
The second derivatives of $\Phi$ and $\Psi$ are bounded thus
\beqs%
&&\hskip-7mm\abs{\RR_{\Phi}\left(\zeta^l_j,\eta^l_j;n^{-\b}
\vvn(s,\frac{j}{n}),n^{-\b}
\uun(s,\frac{j}{n})\right)}+\abs{\RR_{\Psi}\left(\zeta^l_j,\eta^l_j;n^{-\b}
\vvn(s,\frac{j}{n}),n^{-\b} \uun(s,\frac{j}{n})\right)}\\[5pt]
%\\&&
%\phantom{MMMMM}
&&\hskip60mm\leq C \left((\zeta^l_j-n^{-\b}
\uun(s,\frac{j}{n}))^2+
(\eta^l_j-n^{-\b} \vvn(s,\frac{j}{n}))^2 \right),%\\&&
%\abs{\RR_{\Psi}\left(\zeta^l_j,\eta^l_j;n^{-\b}
%\vvn(s,\frac{j}{n}),n^{-\b} \uun(s,\frac{j}{n})\right)}
%%\\ &&
%%\phantom{MMMMM}
%\leq C \left((\zeta^l_j-n^{-\b}
%\uun(s,\frac{j}{n}))^2+ (\eta^l_j-n^{-\b} \vvn(s,\frac{j}{n}))^2
%\right),
\eeqs%
which means that it is sufficient to estimate%
\beqs%
n^{2\b-1}\int_{\Omn} \sum_{j\in \Tn} \left(\zeta^l_j-n^{-\b}
\uun(t,\frac{j}{n})\right)^2 d\mun_t\, \textup{ and }\,
n^{2\b-1}\int_{\Omn} \sum_{j\in \Tn} \left(\eta^l_j-n^{-\b}
\vvn(t,\frac{j}{n})\right)^2 d\mun_t
\eeqs%
uniformly in $t$. We estimate the first expression, the other
will follow the same way. We
denote%
\beqs%
\wzeta_j:=\wzeta_j(t,\uo)=\frac1{l}\sum_{i=0}^{l-1}
\left(\zeta_j-n^{-\b}\uun(t,\frac{j}{n})\right).
\eeqs%
Since $\px \uun(t,x)$ is uniformly bounded for $(t,x)\in
[0,T]\times \T$, we have
\[
(\wzeta_j)^2-(\zeta_j^l-n^{-\b} \uun_j)^2=\Ordo(n^{-\b-1} l)\]
uniformly in $j\in \Tn$, $t\in[0,T]$ and it is enough to estimate %
\beqs n^{2\b-1}\int_{\Omn} \sum_{j\in \Tn}\left(\wzeta_j\right)^2 d\mun_t. \eeqs%
Applying the entropy inequality with respect to the time-dependent
reference measure $\nuun_t$ and using Hölder's inequality:%
\beq%
&&\hskip-10mmn^{2\b-1}\int_{\Omn} \sum_{j\in \Tn}\left(\wzeta_j\right)^2 d\mun_t\leq%\notag \\
%&& \hskip10mm
\frac1{\g} \,n^{2\b-1} H(\mun_t|\nuun_t)+\frac1{\g}\,
 l\,
{n^{2\b-1}} \sum_{j \in \Tn} \log \expect_{\nuun_t} \exp \left(\g
l (\wzeta_j)^2\right),\label{eq:hyp_bl_rpl_1}
\eeq%
for any $\g>0$. $\DD$ is compact, $\zeta$ is bounded thus there
exists a positive constant $C$ such that
\[
\log \expect_{u,v} \exp\big((\zeta-u)y\big)\le C y^2
\]
for all $(u,v)\in \cal{D}$ and $y \in \R$. Thus as a consequence
of Lemma \ref{lemma:hyp_kurschak}, there exists a small $\g>0$ for
which%
\[
\frac1{n} \sum_{j \in \Tn} \log \expect_{\nuun_t} \exp \left(\g l
(\wzeta_j)^2\right)<1.
\]
Substituting into (\ref{eq:hyp_bl_rpl_1}):%
\beqs%
n^{2\b-1}\int\limits_{\Omn} \sum_{j\in \Tn}
\left(\zeta^l_j-n^{-\b} \uun(t,\frac{j}{n})\right)^2 d\mun_t<C
n^{2\b-1} H(\mun_t|\nuun_t)+\Ordo(n^{2\b} l^{-1}).
\eeqs%
Collecting all the estimates, from (\ref{eq:hyp_bl_rpl_0}) we get%
\beqs%
n^{2\b-1} \left(H(\mun_t|\nuun_t)-H(\mun_0|\nuun_0)\right)\le C
n^{2\b-1} \int_0^t  H(\mun_s|\nuun_s)ds+
%\notag\\&&+
\Ordo(n^{-\b} \vee n^{\b-1} l \vee n^{-1-\b} l^3\vee n^{2\b}
l^{-1}).
\eeqs%
Choosing $l$ with
\[
n^{2\b}\ll l \ll n^{\frac{1+\b}{3}}
\]
the error term becomes $o(1)$ and the Theorem follows via the
Grönwall inequality. (If we have the logarithmic-Sobolev
condition, and thus a stronger version of Lemma
\ref{lemma:hyp_obe}, then $l$ can be chosen with $ n^{2\b}\ll l
\ll n^{\frac{1+\b}{2}} $ to make all the error terms $o(1)$.)

\noindent The proof of Lemma \ref{lemma:hyp_kurschak} can be found
in \cite{tothvalko1} or \cite{tothvalko3}.
\begin{lemma}\label{lemma:hyp_kurschak} Suppose
$\xi_1,\xi_2,\dots$ are independent random variables with $\expect
\xi_i=0$ for which
\[
\log \expect \exp (y \xi_i)\le C y^2
\]
with a positive constant $C$ independent of $i$ and $y$. Then
there exists a small positive constants $\g$ depending only on $C$
such that
\[
\log \expect \exp\left(\g l ( \xi^l_1)^2\right)<1.
\]
\end{lemma}

\noindent{\bf\large Acknowledgements:} The author thanks Bálint
Tóth for introducing him into the theory of hydrodynamic limits,
and particularly to the problem of the present paper. He also
thanks Fraydoun Rezakhanlou for drawing his attention to the paper
\cite{dipernamajda} and the method of geometric optics.
\\
This work was partially supported by the Hungarian Scientific
Research Fund (OTKA) grant no. T037685.

%%%%%%%%%%%%%%%%%%%%%%%%%%%%%%%%%%%%%%%%%%%%%%%%%%%%%%%%%%%%%%%%

%\vskip1cm %%
%
%\noindent \small{\textsc{Alfréd Rényi Institute of Mathematics,
%Hungarian Academy of Sciences,\\ Reáltanoda u.~13-15, H-1053
%Budapest, Hungary,\\[5pt]
%Institute of Mathematics, Technical University Budapest,\\ Egry
%J\'ozsef u.~1, H-1111 Budapest, Hungary\\[5pt]
%email: {\tt valko{@}renyi.hu} }}

%%H-1111 Budapest, Hungary\\
%{\tt valko{@}renyi.hu} }}}
%
%%%
%\hbox{\sc
%%\vbox{\noindent
%%\hsize66mm
%%B\'alint T\'oth\\
%%Institute of Mathematics\\
%%Technical University Budapest\\
%%Egry J\'ozsef u. 1.\\
%%H-1111 Budapest, Hungary\\
%%{\tt balint{@}math.bme.hu}
%%}
%%\hskip5mm
%\vbox{\noindent \hsize90mm
%Benedek Valk\'o\\
%Alfréd Rényi Institute of Mathematics\\
%Hungarian Academy of Sciences\\
%Reáltanoda u. 13-15.\\
%H-1053 Budapest, Hungary\\
%%H-1111 Budapest, Hungary\\
%{\tt valko{@}renyi.hu} }}
%

\end{document}